\documentclass{amsart}
\usepackage{amsmath}
\usepackage{amssymb}
\usepackage{amsopn}
\usepackage{amscd}
\usepackage{epsfig}
\usepackage{amsfonts}
\usepackage{latexsym}
\usepackage{graphicx}

\usepackage{pdfsync}

\newcommand{\A}{{\mathcal A}}

\newcommand{\thm}{\tfrac12|\A_m|}
\newcommand{\Pref}{\mathrm{ Pref}}
\newcommand{\Q}{{\mathcal Q}}
\newcommand{\E}{{\mathcal E}}
\newcommand{\Od}{{\mathcal O}}

\theoremstyle{plain}
\newtheorem{theorem}{Theorem}[section]

\newtheorem{lemma}{Lemma}[section]
\newtheorem{proposition}{Proposition}[section]

\title[]{On the structure  of Thue-Morse subwords, with an application to  dynamical systems}

\author[]{F.~Michel Dekking}

\address{ Delft University of Technology.}

\date{\today}

\begin{document}
\maketitle

\begin{abstract} We give an in depth analysis of the subwords of the Thue-Morse sequence. This allows us to prove that
there are infinitely many injective primitive  substitutions with Perron-Frobenius eigenvalue 2 that generate  a
symbolic dynamical system topologically conjugate to the Thue-Morse dynamical system.

\smallskip


\medskip


\end{abstract}

\section{Introduction}

We consider the Thue-Morse sequence $x=0110100110010110\dots$ fixed point of the substitution $\theta$ given by\\[-.6cm]
$$\theta(0)=01,\;\theta(1)=10.$$
By taking its orbit closure under the shift map, the sequence $x$ generates a dynamical system called the Thue-Morse dynamical system.
In the recent paper \cite{CDK2013} it is proved that there are 12 primitive injective substitutions of length 2 that generate
a system topologically conjugate to the Thue-Morse system. A natural question is: what is the list of all primitive injective substitutions
whose incidence matrix has maximal eigenvalue 2 that generate a system topologically conjugate to the Thue-Morse system?

The usual way to generate systems topologically conjugate to a given substitution dynamical system is to consider
 the $N$-block substitution associated to the substitution (\cite{HP}, \cite{CDK2013}).
 See Section~\ref{sec: N-block} for more details, here we give an example: the 5-block substitution $\theta_5$ associated to
 the Thue-Morse substitution $\theta$.

There are twelve Thue-Morse subwords of length $N=5$ (see Example 1 in Section~\ref{sec: comb} for the complete list):
 $w_1=00101,\dots,w_4=01011,\dots,w_{12}=11010$.

 The $\theta_5$-image of a $w_i$ is obtained as the prefix of length 5 of $\theta(w_i)$ followed by
 the prefix of length 5 of $\theta(w_i)$ with the first letter discarded.
For example, since $\theta(00101)=0101100110$, we have  $\theta_5(w_1)=w_4w_{10}$, since $w_{10}=10110$. 
In this way one obtains
\begin{align*}
\theta_5(w_1)&=w_4w_{10},\; \theta_5(w_4)=w_5w_{11},\;\theta_5(w_7)=w_7w_1,\;\theta_5(w_{10})=w_8w_2,\\
\theta_5(w_2)&=w_4w_{10},\;\theta_5(w_5)=w_6w_{12},\; \theta_5(w_8)=w_7w_1,\; \theta_5(w_{11})=w_9w_3,\\
\theta_5(w_3)&=w_5w_{11},\;\theta_5(w_6)=w_6w_{12},\; \theta_5(w_9)=w_8w_2,\; \theta_5(w_{12})=w_9w_3
\end{align*}
We go from this substitution, which is not injective, to an injective one by redistributing the four letters
in the $\theta_5$-images of words of length 2 with odd indices---which always occur in pairs, i.e., the couples $ w_5w_{11},\,w_7w_1$, and $w_9w_3.$
 Concretely, we define a new substitution $\zeta_5$ by keeping $\zeta_5(w_i)=\theta_5(w_i)$ for all words
with an even index, and changing the six others in pairs as, e.g.,
$$\theta_5(w_7)\theta_5(w_1)=w_7w_1\;\,w_4w_{10}=w_7w_1w_4\;\,w_{10}=\zeta_5(w_7)\zeta_5(w_1).$$
This leads to the substitution given by
\begin{align*}
\zeta_5(w_1)&=w_{10},\phantom{w_1}   \;\zeta_5(w_4)=w_5w_{11},\phantom{w_1}   \;\zeta_5(w_7)=w_7w_1w_4, \;\zeta_5(w_{10})=w_8w_2,\\
\zeta_5(w_2)&=w_4w_{10},             \;\zeta_5(w_5)=w_6w_{12}w_9,             \;\zeta_5(w_8)=w_7w_1, \phantom{w_1}  \; \zeta_5(w_{11})=w_3,\\
\zeta_5(w_3)&=w_{11},\phantom{w_1}   \;\zeta_5(w_6)=w_6w_{12},\phantom{w_1}   \;\zeta_5(w_9)=w_8w_2w_5, \; \zeta_5(w_{12})=w_9w_3
\end{align*}
Obviously the substitution $\zeta_5$ is injective, and it is not hard to see that $\zeta^n_5(w_6)=\theta^n_5(w_6)$ for all $n\ge 1$.
Thus, \emph{if}  $\zeta_5$ would be  a primitive substitution, then $\zeta_5$ would generate the same dynamical system as $\theta_5$.
However, $\zeta_5$ is \emph{not} primitive, since  $\zeta_5^2(w_3)=\zeta_5(w_{11})=w_3$.

In Section \ref{sec: injective} we will repair this defect by defining a substitution $\eta_5$ which generates
 the same dynamical system as $\theta_5$, but \emph{is} primitive. Actually, we give this construction for all $\eta_N$,
 where $N$ is a power of two plus one. For this we need an explicit expression for $\theta_N$, which is given in Section \ref{sec: N-block},
 based on the combinatorial analysis in Section \ref{sec: comb}. Our main result is in Section \ref{sec: morse}: there exist
 infinitely many substitutions in the Thue-Morse conjugacy class if we allow also
 non-constant length substitutions with Perron-Frobenius eigenvalue 2.

\section{Combinatorics of Thue-Morse subwords}\label{sec: comb}

The subwords of the Thue-Morse sequence have been well studied (see, e.g., \cite{ACS1997}).
We  show here  that the subwords of length $N=2^m+1$ have a particularly elegant structure for $m=2,3,\dots$.
Let $\A_m$ be the set of these words. It is well known (and will be reproved here) that the cardinality of $\A_m$ equals $|\A_m|=3\cdot2^m$.
 We lexicographically order the words in $\A_m$, representing them as $$w^m_1<w^m_2<\dots<w^m_{|\A_m|}.$$
Crucial to the following analysis is the partition of $\A_m$ into 4 sets
$$\A_m=\Q_1\cup\Q_2\cup\Q_3\cup\Q_4,$$
where each  $\Q_k$ consists of one quarter of consecutive words from $\A_m$.
If we want to emphasize the dependence on $m$ we write $\Q_k^m$. Let
$$q_k= \min \Q_k, \quad\mathrm{for}\; k=1,2,3,4.$$
Thus\\[-.6cm]
 $$q^m_1=w^m_1,\quad q^m_2=w^m_{\frac14|\A_m|+1},\quad q^m_3=w^m_{\frac12|\A_m|+1},\quad q^m_4=w^m_{\frac34|\A_m|+1}.$$

Let $f_0^\omega=0110\dots$ and $f_1^\omega=1001\dots$ be the two infinite fixed points of $\theta$, and let $f_0=f_0^m$ and
$f_1=f_1^m$ be the length $2^m+1$ prefixes of $f_0^\omega$ and $f_1^\omega$.

\medskip

\noindent {\bf Example 1}\; The case $m=2$. The set $\A_2$ is given by\\
$\{00101,00110,01001,01011,01100,01101,10010,10011,10100,10110,11001,11010\}.$\\
Here $q_1=00101, q_2=01011, q_3=10010, q_4=10110,$ and $f_0=01101,f_1=10010$.

\medskip

We use frequently mirror invariance of the Thue-Morse words, i.e., if the mirroring operation is define as the length 1
substitution given by $\widetilde{0}=1, \widetilde{1}=0$, then $u$ is a Thue-Morse subword if and only if $\widetilde{u}$ is
a Thue-Morse subword. This follows directly from $\widetilde{\theta(0)}=\theta(1)$.

The Thue-Morse substitution $\theta$ has the following trivial, but important property.

\begin{lemma}\label{lem: order} If words $u$ and $v$ satisfy $u<v$, then $\theta(u)<\theta(v).$
\end{lemma}

The words in $\A_{m+1}$ are generated by the words in $\A_m$ in a simple way.
Each word $u\in\A_m$ has two $\delta\varepsilon$scendants, $\delta(u)$ and $\varepsilon(u)$,
where, by definition, $\delta(u)$ is the length $2^{m+1}\!+1$ prefix, and $\varepsilon(u)$ the length $2^{m+1}\!+1$ sufffix
of $\theta(u)$. For example: since $\theta(00101)=0101100110$, we have
$$\delta(00101)=010110011\in\A_3,\quad \varepsilon(00101)=101100110\in\A_3.$$
The next lemma follows from Lemma \ref{lem: order}.
\begin{lemma}\label{lem: descorder} If two words $u$ and $v$ satisfy $u<v$, then $\delta(u)<\delta(v)$. If moreover, $u_1=v_1$, then $\varepsilon(u)<\varepsilon(v)$.
\end{lemma}

In the following we will freely use group notation for words over the alphabet $\{0,1\}$. For instance $(01)^{-1}0110=10$.

\begin{proposition}\label{prop: qandf} For all $m$ the smallest words in the $\Q_k^m$ can be expressed in $f_1^m$:
$$(1)\; q_1=1^{-1}f_1\,1,\quad  (2)\; q_2=(10)^{-1}f_1\,11,\quad (3)\; q_3=f_1,\quad (4)\; q_4=(100)^{-1}f_1\,110.$$
\end{proposition}

This proposition is tied up with the following one.

\begin{proposition}\label{prop: Qprop} For all $m=2,3,\dots$
\begin{eqnarray*}
(1)& \Q_1^{m+1}=\varepsilon(\Q_3^m\cup\Q_4^m), \quad (2)\quad \Q_2^{m+1}=\delta(\Q_1^m\cup\Q_2^m),\\
(3)& \Q_3^{m+1}=\delta(\Q_3^m\cup\Q_4^m),    \quad (4)\quad \Q_4^{m+1}=\varepsilon(\Q_1^m\cup\Q_2^m).
\end{eqnarray*}
\end{proposition}

We first prove (3) of Proposition~\ref{prop: qandf}. By mirror symmetry of the Thue-Morse words we know that exactly half of the
words in $\A_m$ start with 1, so $q_3^m$ is the \emph{smallest} word starting with 1 in $\A_m$.  From the example above we see that $q_3^2=f_1^2$.
Then it  follows by induction  that
$f_1^m$ is \emph{also} the smallest word  with prefix 1, since $\delta$ is order preserving,
and the other words starting with 1 in $\A_{m+1}$ are generated by words with prefix 00 or 01,
which under $\varepsilon$ generate words with prefix 101 or 110, which are both larger than the prefix 100 from $f_1^{m+1}$.

\medskip

\emph{Proof of Proposition~\ref{prop: Qprop}:}
We first prove (3). By Proposition~\ref{prop: qandf} (3) the smallest symbol in $\Q_3^{m+1}$ is mapped by $\delta$ to the smallest symbol
 of $\delta(\Q_3^m\cup\Q_4^m)$. But since $\delta$ is orderpreserving, it follows by matching cardinalities that (3) holds.

Since $\Q_1^m\cup\Q_2^m$ maps under $\delta$ to words starting with 0, where the largest word is $\delta(f^m_0)=f^{m+1}_0$,
again a cardinality argument shows that its image must be $\Q_2^{m+1}$, so (2) holds.

Since $\varepsilon$ maps consecutive symbols starting with 0 to consecutive symbols starting with 1,
$\varepsilon(\Q_1^m\cup\Q_2^m)$ must be $\Q_4^{m+1}$.
Then for $\varepsilon(\Q_3^m\cup\Q_4^m)$ there is only $\Q_1^{m+1}$ left, i.e., (1) holds.\qed

\medskip

\emph{Proof of Proposition~\ref{prop: qandf}:}
We already proved (3), i.e., that $q_3^m=f_1^m$. From this it follows that $1^{-1}f_1^m$ is smaller than (or equal to) all words of length $2^m$ starting with 0, except maybe those that are not of the form $1^{-1}w$ with $w\in\A_m$. But these are of the form $0^{-1}w$, where $w$ starts with 00. This implies that $0^{-1}w$ starts with $01>00$, since 000 does not occur in a Thue-Morse word. Conclusion: $1^{-1}f_1^m$ is the smallest of all words \emph{of length} $2^m$ starting with 0. It has a unique right extension to the word $1^{-1}f_1^m\,1$, which is still the smallest among all words of length $2^m+1$, i.e., (1) holds.

To prove (2), note that $q_2^{m+1}=\delta(q_1^m)$ by (2) of Proposition~\ref{prop: Qprop}. Also note that $\theta(q_1^m)$ has suffix 0, and $\theta(f_1^m)$ has suffix 1 for all $m$. Applying $\theta$ to both sides of (1) we obtain
$$q_2^{m+1}=\theta(q_1^m)0^{-1}=(10)^{-1}\theta(f_1^m)10\,0^{-1}=(10)^{-1}f_1^{m+1}11.   $$

To prove (4), note that $q_4^{m+1}=\varepsilon(q_1^m)$, by  Proposition~\ref{prop: Qprop} (4). It follows that
$$q_4^{m+1}=0^{-1}\theta(q_1^m)=0^{-1}\theta(1^{-1}f_1^m\,1)=0^{-1}(10)^{-1}\theta(f_1^m)10=(100)^{-1}f_1^{m+1}110.\qed $$

\medskip

We would like to make the following historical remarks. Our Proposition~\ref{prop: qandf} (3) is the finite, mirrored,  version of
 Corollaire 4.4 in \cite{Berst1994} by Berstel.
Our Proposition~\ref{prop: qandf} (1) is the finite  version of
 Corollary 2 to Theorem 1 in \cite{ACS1997} by Allouche, Curry and Shallitt.

\medskip

\noindent In the next section we will need the following lemma, in which we use some new notation. For a word $w=w_1\dots w_k$, we write $\Pref_{\ell}(w)$ for
its prefix $w_1\dots w_\ell$ of length $\ell\le k$.

\begin{lemma}\label{lem: firsthalf} For all $m\ge 1$ and $N=2^m+1$ we have
$$\Pref_N(w^{m+1}_{2i-1})=\Pref_N(w^{m+1}_{2i})=w^m_i, \quad \mathrm{for}\;i=1,\dots,|\A_m|.$$
\end{lemma}

\emph{Proof:} Note first that \emph{all} words $w^m_i$ have to appear as an $N$-prefix of the words $w^{m+1}_{j}$, and in \emph{lexicographical order}.
Here $\le$ can, and will occur, and the only fact that has to be checked is that there are no words $w^m_i$ occurring  only once.

 A quick glance at $\A_2$ in the example above  shows this is true for $m=1$, since  $\A_1$ is equal to $\{001,010,011,100,101,110\}$.
Suppose it is true for $m$. Then for all $i$ the two words $w^{m+1}_{2i-1}$ and $w^{m+1}_{2i}$ will
have two $\delta$-descendants that have the same prefix of length $N+1$, and the same holds for the two $\varepsilon$-descendants.
So there are no words $w^{m+1}_j$ occurring  only once as a $(N+1)$-prefix of a word $w^{m+2}_k$.
 \qed

 \medskip

\noindent {\bf Example 2}\; The case $m=3$. The set $\A_3$ has 24 elements  given by\\
\begin{tabular}{rrrr} $w_{1} =0 0 1 0 1 1 0 0 1$ &        $ w_{7} =0 1 0 1 1 0 0 1 1   $ &   $w_{13} = 1 0 0 1 0 1 1 0 0$  &      $w_{19} =1 0 1 1 0 0 1 1 0$\\
$w_{2} =0 0 1 0 1 1 0 1 0$  &       $ w_{8} =0 1 0 1 1 0 1 0 0   $ &   $w_{14 } =1 0 0 1 0 1 1 0 1$  &      $w_{20} =1 0 1 1 0 1 0 0 1$\\
$w_{3} =0 0 1 1 0 0 1 0 1$  &       $ w_{9} =0 1 1 0 0 1 0 1 1 $ &     $w_{15 } =1 0 0 1 1 0 0 1 0$  &      $w_{21} =1 1 0 0 1 0 1 1 0$\\
$w_{4} =0 0 1 1 0 1 0 0 1$  &       $w_{10} =0 1 1 0 0 1 1 0 1$ &      $w_{16 } =1 0 0 1 1 0 1 0 0$  &      $w_{22} =1 1 0 0 1 1 0 1 0$\\
$w_{5} =0 1 0 0 1 0 1 1 0$  &       $w_{11} =0 1 1 0 1 0 0 1 0$ &      $w_{17 } =1 0 1 0 0 1 0 1 1$  &      $w_{23} =1 1 0 1 0 0 1 0 1$\\
$w_{6} =0 1 0 0 1 1 0 0 1$  &       $w_{12} =0 1 1 0 1 0 0 1 1$ &      $w_{18 } =1 0 1 0 0 1 1 0 0$  &      $w_{24} =1 1 0 1 0 0 1 1 0.$\\
\end{tabular}

\noindent  Here $q_1=w_1, q_2=w_7, q_3=w_{13}, q_4=w_{19},$ and $f_0=w_{12},f_1=w_{13}$.
\qed

\section{The Thue-Morse $N$-block substitutions $\theta_N$ }\label{sec: N-block}

 A simple way to produce substitutions that generate dynamical systems topologically conjugate to a given substitution
 is to construct $N$-block substitutions---see Section 4 of \cite{CDK2013}.

 We will describe this construction for a general substitution $\alpha$ of constant length  on an alphabet $A$. Let
the length of $\alpha$ be $L$, an integer greater  than one. Further, let $N$ denote any positive integer.
Let $ {\mathcal L}_{\alpha}$ be language of $\alpha$, i.e., the collection of all words occurring in some power $\alpha^n(a)$, for some $a\in A$.
We define the alphabet $B=A^N\cap\,{\mathcal L}_{\alpha}$, and construct a substitution
$\alpha_N$ on $B$, called the \emph{$N$-block substitution} associated to $\alpha$. Namely, if $b=a_1\dots a_N$ is an element of $B$, we apply $\alpha$ to $b$, obtaining a word $v:=\alpha(a_1\dots a_N)$ of length $LN$.
We then define\\[-.6cm]
$$\alpha_N(b)=v_1\dots v_N, v_2\dots v_{N+1},\dots,v_{L}\dots v_{L+N-1}.$$

\noindent {\bf Example 3} Let $N=3$,  let $A=\{0,1\}$, and let $\alpha=\theta$, the Thue Morse substitution.
Then the words of length $N$ in the language of $\theta$ are $w_1=001,\dots,w_6=110$.
Since $\theta(001)=010110$, we have $\theta_3(w_1)=w_2w_5$, and similarly we find \\[-0.6cm]
\begin{align*}
\theta_3(w_1)&=w_2w_5,\; \theta_3(w_2)=w_3w_6,\;\theta_3(w_3)=w_3w_6,\\ \theta_3(w_4)&=w_4w_1,\;\theta_3(w_5)=w_4w_1,\; \theta_3(w_6)=w_5w_2.\qed
\end{align*}

\smallskip

\noindent We give an explicit formula for all $\theta_N$, where $N=2^m+1$ for $m=2,3,\dots.$

It is convenient to define the translation $\tau$ on $\{1,\dots,|\A_m|\}$ by
$$\tau(i)=\big(i-1+\thm\big)\!\!\!\mod|\A_m|\,+\,1.$$
We extend $\tau$ to $\tau:\A_m\rightarrow\A_m$ by putting $\tau(w_i)=w_{\tau(i)}$.

\begin{proposition}\label{prop: Nblock} Let $\theta$ be the Thue-Morse substitution, and let $N=2^m+1$, with $m\ge2$. Write $\theta_N(w_i)=w_{F(i)}w_{G(i)}$ for $i\in \A_m$. Then
\begin{align}
F(2i)&=F(2i-1)=\tfrac14|A_m|+i& \qquad\mathrm{for}\; i=1,\dots,\thm.\\
G(i)&=\tau(F(i))& \qquad\mathrm{for}\; i=1,\dots,|\A_m|.
\end{align}
\end{proposition}

\emph{Proof:} Note that $w_{F(i)}=\Pref_N[\delta(w_i)]$, and $w_{G(i)}=\Pref_N[\varepsilon(w_i)]$.

We first show that $\theta_N(q_1)=q_2q_4$. This follows directly from Proposition~\ref{prop: Qprop} (2) and (4), since by Proposition~\ref{prop: qandf} (2) and (4) we have $\Pref_N[ q_2^{m+1}]=q_2^m$, and $\Pref_N[ q_4^{m+1}]=q_4^m$.
 So (1) holds for $i=1$. Similarly, we have $\theta_N(q_3)=q_3q_1$.

It follows directly from Lemma~\ref{lem: firsthalf} that $F(2i)=F(2i-1)$ for all $i$, and since $\delta$ is orderpreserving (1) follows from the $i=1$ case.

In the same way (2) follows from the $i=1$ and the $i=\thm+1$ case. \qed

\section{the construction of injective Thue Morse substitutions}\label{sec: injective}

 The substitution $\theta_N$ is exactly 2-to-1. In this section we construct for  $m\ge 3$ a 1-to-1
 substitution $\eta_N$ on $\A_m$ which admits one of the fixed points of $\theta_N$ as a fixed point.
 The idea for this construction is a sort of converse of a construction in \cite{Dek78}.

 \noindent Notationally it is convenient to introduce the set $\E_m$ of words with even indices, and the set $\Od_m$ of words with odd indices.

 The substitution $\eta_N$  will be a non-constant length substitution with lengths 1, 2 or 3.
 It is defined by $\eta_N(w_{i})=\theta_N(w_{i})$ for $w_i\in \E_m$, and
\begin{equation}\label{eq: etaN}
 \eta_N(w_{i})=
\begin{cases} w_{G(i)}        &\quad\mathrm{for}\; w_i\in \Od_m\cap\Q_1,\\
w_{F(i)}                      &\quad\mathrm{for}\; w_i\in \Od_m\cap\Q_2,\\
\theta_N(w_{i})w_{F(\tau( i))}  &\quad\mathrm{for}\; w_i\in \Od_m\cap\Q_3,\\
w_{G(\tau( i))}\theta_N(w_{i})  &\quad\mathrm{for}\; w_i\in \Od_m\cap\Q_4.
\end{cases}
\end{equation}

The  idea of this definition is that $\theta_N$ and $\eta_N$ act in the same way on words of length 2 occurring at even places in the fixed point $f_0^\omega$ of $\theta_N$. Suppose for instance that $w_i\in \Od_m\cap\Q_2$.
Then by (1) of Proposition \ref{prop: Nblock} there is a unique $w_j\in \Od_m\cap(\Q_1\cup\Q_2)$ such that $F(j)=i$.
Note that   $w_{G(j)}\in \Od_m\cap\Q_4$, since  by (2) of Proposition \ref{prop: Nblock},  $G(j)=\tau(F(j))=\tau(i)$,
and $\tau(\Q_2)=\Q_4$.
Therefore for all odd $j$ with $w_j\in  \Q_1 \cup \Q_2$  and since $\tau$ is an involution,
\begin{align*} \eta_N(w_{F(j)}w_{G(j)})&= \eta_N(w_{i}w_{\tau(i)})=w_{F(i)}\,w_{G(i)}\theta_N(w_{\tau(i)})\\
&=  \theta_N(w_{i}w_{G(j)})=   \theta_N(w_{F(j)}w_{G(j)}).
\end{align*}
Similarly, if $w_i\in \Od_m\cap\Q_3$, then there is a unique $w_j\in \Od_m\cap(\Q_3\cup\Q_4)$ such that $F(j)=i$.
Now $w_{G(j)}\in \Od_m\cap\Q_1$, and we have for all odd $j$ with $w_j\in  \Q_3 \cup \Q_4$
\begin{align*} \eta_N(w_{F(j)}w_{G(j)})&= \eta_N(w_{i}w_{\tau(i)})= \theta_N(w_{i})w_{F(\tau(i))}\,w_{G(\tau(i))}\\
&=\theta_N(w_{i}w_{\tau(i)})=   \theta_N(w_{F(j)}w_{G(j)}).
\end{align*}

\noindent Since $\eta_N(w_{i})=\theta_N(w_{i})$ for $w_i\in \E_m$, it follows that for all $j\in\A_m$
$$\eta_N(w_{F(j)}w_{G(j)})=\theta_N(w_{F(j)}w_{G(j)}).$$
\noindent Since $w_{F(j)}$ is \emph{always} followed by $w_{G(j)}$, it must be that
$$\eta_N^n(f^m_0)=\theta_N^n(f^m_0)\quad \mathrm{for\;} n=1,2,\dots.$$
If we knew that $\eta_N$ was primitive, i.e., its incidence matrix is primitive, then this would imply that
 $\eta_N$ and $\theta_N$ generate the same minimal set: $X_{\eta_N}=X_{\theta_N}.$

\section{The problem of primitivity}\label{sec: prim}

\begin{proposition}\label{prop: prim} The substitution $\eta_N$ is primitive.
\end{proposition}

\emph{Proof:} Let us write $v\twoheadrightarrow w$ for $v,w\in\A_m$ if there exists an $n$ such that $w$ occurs in $\eta_N^n(v)$.
We will prove the following:
$$ \mathrm{(F)} \; f_0\twoheadrightarrow w, \;f_1\twoheadrightarrow w \mathrm{\;\;for\; all\;} w\in\A_m,\;\;
  \mathrm{ (B)}  \;\mathrm{ either}\; v\twoheadrightarrow f_0, \;\mathrm{or}\;v\twoheadrightarrow f_1 \mathrm{\;\;for\; all\;} v\in\A_m.$$
Obviously (B)+(F) implies that $v\twoheadrightarrow w$ for all $v,w\in\A_m$, i.e., $\eta_N$ is irreducible.
But primitivity follows easily from this, by observing that $f_0\twoheadrightarrow f_0$ and $f_1\twoheadrightarrow f_1$
in \emph{one} step.

\noindent For the proof of (F), note that $\eta_N^n(f^m_0)=\theta_N^n(f^m_0)$ implies that every $w$ occurs in some $\eta_N^n(f^m_0)$,
since $\theta_N$ is primitive. For $\eta_N^n(f^m_1)$ such an equality does not hold, but it is still true that $\eta_N^n(f^m_1)$ is
a prefix of $\theta_N^n(f^m_1)$, which of course leads to the same conclusion.

\noindent The proof of (B) is somewhat more involved, and we will prove something stronger, namely that starting from any $v$ either
$f_0$ or $f_1$ will occur as \emph{first} letter of $\eta_N^n(v)$ for some $n\ge1$.
We first study $\theta_N$,  defining the `initials' map $\phi:\A_m\rightarrow\A_m$ by
$$\phi(w_i)=w_{F(i)}, \;\;\mathrm{where}\;\;w_{F(i)}=\Pref_1(\theta_N(w_i)) \quad \mathrm{for}\; i=1,2,\dots,|\A_m|.$$
From Proposition \ref{prop: Nblock}    we obtain
$$\phi(\Q_1)\subseteq\phi(\Q_2),\quad\phi(\Q_2)\subseteq\phi(\Q_2),\quad\phi(\Q_3)\subseteq\phi(\Q_3),\quad\phi(\Q_4)\subseteq\phi(\Q_3).$$
Moreover,  $\phi$ is strictly increasing on $\Q_2\setminus\{f_0\}$, which implies that
$$\phi^n(w_i)=f_0 \quad\mathrm{for\; all}\;n\ge\thm, \;w_i\in \Q_1\cup\Q_2.$$
By mirroring, we have $\phi^n(w_i)=f_1$ for all $n\ge\thm$ and $w_i\in \Q_3\cup\Q_4.$

Next, we define $\psi:\A_m\rightarrow\A_m$ by
$$\psi(w_i)=w_{H(i)}, \;\;\mathrm{where}\;\;w_{H(i)}=\Pref_1(\eta_N(w_i)) \quad \mathrm{for}\; i=1,2,\dots,|\A_m|.$$
From the definition of $\eta_N$ we see that $\psi|_{\Q_2\cup\Q_3}=\phi|_{\Q_2\cup\Q_3}$, which implies that $\psi^n(\Q_2\cup\Q_3)\subseteq \{f_0,f_1\}$ for $n\ge\thm$.

Note that  $\psi(\Q_1)\subseteq \Q_4$, so what remains is to study the behavior of $\psi$ on $\Q_4$.
First, if $w_i\in \Q_4\cap\E_m$, then $\psi(w_i)=\phi(w_i)\in\Q_3$, so $\psi^n(w_i)=f_1$ for all large $n$.

Second, if $w_i\in \Q_4\cap\Od_m$, then $\psi(w_i)>w_i$. So  iterating $\psi$ will always result in hitting a letter $w_j$ with
$w_j\in \Q_4\cap\E_m$, and we are in the first case. \qed

 \section{An infinite Thue-Morse conjugacy list}\label{sec: morse}

 Our work in the previous sections leads to an answer to the Thue-Morse conjugacy  list question.

\begin{theorem}
There are infinitely many injective primitive substitutions with Perron-Frobenius eigenvalue 2 that generate a dynamical system
topologically conjugate to the Thue-Morse dynamical system.
\end{theorem}

\emph{Proof:} Infinitely many substitutions are given by the $\eta_N$, where $N=2^m+1$ for $m=2,3,\dots.$
These generate minimal systems conjugate to the Thue-Morse substitutions
because they generate the same systems as the $N$-block substitutions $\theta_N$.
 From Proposition \ref{prop: Nblock} and the defining Equation (\ref{eq: etaN}) it follows that the $\eta_N$ are injective, and primitivity is given by
 Proposition \ref{prop: prim}.

It remains to prove that the Perron-Frobenius eigenvalue of the matrix $M_N$ of $\eta_N$ is equal to 2.
This is clear from the definition, but here is a formal proof. Let $e_N$ be the vector of length $|\A_m|$ with all ones, and let $d_N$ be the
vector of length $|\A_m|$  with all zero's, except for a 1 at the position of $f^m_0$. Let $\ell(v)$ be  the length of a word $v$.
Then for all $n\ge 1$ one has $ d_N^{\sc{T}}\,M^n_N\,e_N=\ell(\eta^n_N(f^m_0))=\ell(\theta^n_N(f^m_0))=2^n.$
Since $M_N$ is a primitive non-negative matrix it follows from the Perron-Frobenius theorem
 that this implies that the eigenvalue of largest modulus is equal to 2. \qed

\bibliographystyle{plain}

 \end{document}